\documentclass[12pt]{article}
\title{Categoricity for  Patterns of Order 2}
\author{Timothy J. Carlson}

\date{}

\usepackage{amsmath}    
\usepackage{amssymb}
\usepackage{graphicx}   
\usepackage{verbatim}   
\usepackage{color}      
\usepackage{subfigure}  
\usepackage{hyperref}   

\newcommand{\qed}{\hfill {\bf QED}}
\newcommand{\blankline}{\vspace{4 mm}}
\newcommand{\halfblankline}{\vspace{2 mm}}

\newcommand{\bfP}{\mbox{${\bf P}$}}
\newcommand{\bfPi}{\mbox{${\bf P}_i$}}

\newcommand{\bfPiplusone}{\mbox{${\bf P}_{i+1}$}}
\newcommand{\bfPinfty}{\mbox{${\bf P}_\infty$}}

\newcommand{\bfPn}{\mbox{${\bf P}_n$}}

\newcommand{\bfPnstar}{\mbox{${\bf P}_n^*$}}
\newcommand{\bfPnplusone}{\mbox{${\bf P}_{n+1}$}}

\newcommand{\bfPplus}{\mbox{${\bf P}^+$}}
\newcommand{\bfPprime}{\mbox{${\bf P}'$}}
\newcommand{\bfPstar}{\mbox{${\bf P}^*$}}

\newcommand{\bfPzero}{\mbox{${\bf P}_0$}}

\newcommand{\bfQ}{\mbox{${\bf Q}$}}

\newcommand{\bfQplus}{\mbox{${\bf Q}^+$}}
\newcommand{\bfQprime}{\mbox{${\bf Q}'$}}
\newcommand{\bfQstar}{\mbox{${\bf Q}^*$}}

\newcommand{\calB}{\mbox{${\cal B}$}}

\newcommand{\calL}{\mbox{${\cal L}$}}

\newcommand{\calLtwo}{\mbox{${\cal L}_2$}}

\newcommand{\calR}{\mbox{${\cal R}$}}

\newcommand{\calRone}{\mbox{${\cal R}_1$}}

\newcommand{\calRtwo}{\mbox{${\cal R}_2$}}

\newcommand{\leqone}{\mbox{$\leq_1$}}
\newcommand{\leqpw}{\mbox{$\leq_{pw}$}}
\newcommand{\leqtwo}{\mbox{$\leq_2$}}

\newcommand{\preq}{\mbox{$\preceq$}}

\newcommand{\preqzero}{\mbox{$\preceq_0$}}

\newcommand{\preqone}{\mbox{$\preceq_1$}}

\newcommand{\preqtwo}{\mbox{$\preceq_2$}}

\newcommand{\preqkinfty}{\mbox{$\preceq_k^\infty$}}

\newcommand{\unibfP}{\mbox{$|{\bf P}|$}}

\newcommand{\unibfPplus}{\mbox{$|{\bf P}^+|$}}

\newtheorem{prop}{Proposition}[section]
\newtheorem{cor}[prop]{Corollary}

\newtheorem{lem}[prop]{Lemma}

\newtheorem{dfn}[prop]{Definition}

\begin{document}

\maketitle

\begin{center}
  The Ohio State University, Columbus, OH 43210 USA \\
  email: carlson.6@asc.ohio-state.edu
\end{center}

In this paper we show how a Categoricity Theorem for patterns of resemblance of order 2, in  analogy to Theorem 9.1 of [\ref{refC1}] for \calRone, follows from [\ref{refC2}]. This is the result alluded to in the last paragraph of the introduction to [\ref{refC2}] where it is stated
\begin{itemize}
\item[]
{\it ... a method of generating the core is established which shows that the order in which patterns of embeddings of this level occur is the same for reasonable hierarchies.}
\end{itemize}
As a consequence, if a reasonable hierarchy \calB\ (see the Categoricity Theorem below) has arbitrary long finite chains in the interpretation of \preqtwo\ then a finite structure is a pattern of resemblance of order two iff it is isomorphic to a finite substructure of \calB\   (see Corollary \ref{c characterization}).
These results apply to the version of \calRtwo\ defined in the introduction to [\ref{refC2}] as initial segments are reasonable hierarchies.

Our basic reference is [\ref{refC2}]. 

We will work in the theory ${\sf KP}\omega$ i.e. Kripke-Platek Set Theory plus the Axiom of Infinity.

Fix a language \calL\ including the binary relation symbol \preq. 
Let \calLtwo\ be the expansion of \calL\ by binary relation symbols \preqone\ and \preqtwo.
We also write \preqzero\ for \preq. 
We use {\it structure} to refer to what is more commonly called a partial structure where the interpretations of the function symbols are allowed to be partial. 
We will write $|\bfP|$ for the universe of a structure \bfP.

For the remainder of the paper, let \calR\ be an EM structure (see Section 3 of [\ref{refC2}]) for \calL\  on the class of ordinals with $\preq^{\cal R}$ the usual ordering. 
We assume  the restriction of \calR\ to any ordinal is a set, there is a restriction with $\omega$ indecomposables and the indecomposables are cofinal in the ordinals. 
Since \calR\ is an EM structure, it can be recovered by its restriction to the $\omega^{th}$ indecomposable which implies the set of indecomposables is $\Delta$-definable and the function which maps an indecomposable $\lambda$ to $\calR \upharpoonright \lambda$ is $\Sigma$-definable. 

We also assume \calB\ is a structure for \calLtwo\ whose arithmetic part (i.e. restriction to \calL) is an arithmetic structure with respect to \calR\ (Definition 4.1 of [\ref{refC2}]) in which the interpretation of each function symbol is total. 
We do not require that  \calB\ be well-ordered with respect to the interpretation of \preq\ though our main focus will be on those \calB\ which are.
Recall that $\preq_k^{\cal B}$ respects $\preq_{k-1}^{\cal B}$ if 
$$\alpha\,\preq_{k-1}^{\cal B}\,\beta\, \preq_{k-1}^{\cal B}\, \gamma \ {\rm and}  \  \alpha\, \preq_k^{\cal B}\,  \gamma    \ \ \ \Longrightarrow \ \ \ \alpha\, \preq_k^{\cal B}\,  \beta$$
for all $\alpha, \beta, \gamma$.

\blankline

\noindent{\bf Categoricity Theorem for \calRtwo.}
{\it 
If
\begin{enumerate}
 \item[(a)] 
 For $k=1,2$, $\calB \! \upharpoonright \! \alpha\ \preq_k^\infty \ \calB \! \upharpoonright \! \beta$ whenever 
  $\alpha \preq_k^{\cal B} \beta$.
 \item[(b)] $\preq_1^{\cal B}$ and $\preq_2^{\cal B}$ are partial orderings of the universe of \calB\ with $\preq_2^{\cal B}\subseteq \preq_1^{\cal B}\subseteq \preq_0^{\cal B}$.
 \item[(c)] $\preq_k^{\cal B}$  respects $\preq_{k-1}^{\cal B}$ for $k=1,2$.
 \item[(d)] The arithmetic part of \calB\ is $\calR \! \upharpoonright \! \lambda$ for some $\lambda$ which is indecomposable in \calR.
 \end{enumerate}
then the core of \calB\ is isomorphic to an initial segment of the
core of $\calRtwo \! \upharpoonright \! \lambda$.}

\blankline

\calRtwo\ is defined in Definition 5.4 of [\ref{refC2}].

For the rest of the paper, assume \calB\ satisfies (a)-(c) of the theorem. 
We do not assume that \calB\ is necessarily well-ordered by $\preq^{\cal B}$.

\begin{dfn}
A pattern \bfP\ is {\bf $\cal B$-covered} if there is a covering of \bfP\ in \calB.
\end{dfn}

See Definition 5.3 of [\ref{refC2}] for the definition of {\it pattern} (short for {\it pattern of resemblance of order two}).
See Definition  5.2 of [\ref{refC2}] for the definition of {\it covering}. That definition is slightly different from that used in [\ref{refC1}] in that the range of a covering is required to be closed (Definition 2.3 of [\ref{refC2}]).

\begin{dfn}
\label{d extends above}
Assume ${\bf P}$ is a pattern,
 $h$ is a function from the universe of $\bf P$ into the universe of  \calB\ and $\varphi$ is a regressive
function on the nonminimal
indecomposable elements in the range of $h$ (i.e. $h(\alpha)<\alpha$ for any nonminimal element in the range of $h$ which is indecomposable in \calR). 
Suppose also that  ${\bf P}$ is a closed substructure of the pattern ${\bf P}^+$.  
A function $h^+$ of the universe of ${\bf P}^+$ into the universe of \calB\
{\bf extends $h$ above $\varphi$} if $h^+$ extends $h$ and 
\begin{center}
 $\varphi(h(a))< h^+(b)$
\end{center}
for  any indecomposable $b$ in ${\bf P}^+$  and any indecomposable
$a$ in $\bf P$ such that $(-\infty,a)^{\bf P}\prec^{{\bf P}^+} b \prec^{{\bf P}^+}a$.
\end{dfn}

\begin{dfn}
Assume ${\bf P}$ and ${\bf P}^+$ are patterns and ${\bf P}$ is a closed substructure of 
${\bf P}^+$.
The rule ${\bf P}|{{\bf P}^+}$ is
{\bf cofinally valid} in \calB\ if for every covering $h$ of $\bf P$ in \calB\ and every 
regressive function $\varphi$ on the nonminimal
indecomposable elements in the range of $h$
there is a covering $h^+$ of ${{\bf P}^+}$ into \calB\ which extends
$h$ above $\varphi$.
\end{dfn}

\begin{lem}
Every generating rule is cofinally valid in \calB.
\end{lem}
{\bf Proof.}
The only properties of \calRtwo\ used in the proof  of part 2 of Lemma 13.11 of [\ref{refC2}] and the supporting lemmas are the preliminary properties we have assumed of \calB. Therefore, the proof carries over with \calRtwo\ replaced by \calB. 

The proof is by induction on the generation of the generating rules (Definition 13.10 of [\ref{refC2}]).

Suppose \bfPplus\ is 1-correct arithmetic extension of \bfP. 
The proof that $\bfP | \bfPplus$ is cofinally valid in \calB\ is analogous to  the proof of Lemma 8.4 of [\ref{refC2}].
Assume $h$ is a covering of \bfP\ in \calB\ and $\varphi$ is a regressive function on the nonminimal indecomposables in the range of $h$. 
Notice that any covering of \bfPplus\ in \calB\ which extends $h$ vacuously extends $h$ above $\varphi$ since there are no new indecomposable elements (by Lemma 4.9 of [\ref{refC2}]). 
By Lemma 4.5 of [\ref{refC2}], there is an embedding $h^+$ of the arithmetic part of \bfPplus\ in \calR\ which extends $h$. 
Clearly, the range of $h^+$ is contained in $\lambda$.
A straightforward argument using the fact that \bfPplus\ is a 1-correct arithmetic extension of \bfP\ (Definitions 4.8, 7.1 and 8.1 of [\ref{refC2}]) shows $h^+$ is a covering of \bfPplus\ in \calB.

Suppose \bfPplus\ is obtained from \bfP\ by 1-reflecting $X$ downward from $b$ to $a$.
The proof that $\bfP | \bfPplus$ is cofinally valid in \calB\ is analogous to  the proof of Lemma 9.3 of [\ref{refC2}].
Assume $h$ is a covering of \bfP\ in \calB\ and $\varphi$ is a regressive function on the nonminimal indecomposables in the range of $h$.
Since $h$ is a covering and $a\prec_1^{\bf P} b$, $h(a)\prec_1^{\cal B} h(b)$ implying $h(a)\prec_1^\infty h(b)$ in \calB. 
Therefore, there is $\tilde{X}$ such that 
$h[(-\infty,a)^{\bf P}]\cup \{\varphi(h(a))\} \ < \ \tilde{X} < h(b)$ and $h[(-\infty,a)^{\bf P}]\cup \tilde{X}$ is both closed and a covering of $h[(-\infty,a)^{\bf P}]\cup h[X]$. 
Let $h^+$ be the order isomorphism of \bfPplus\ and $h[|\bfP|]\cup \tilde{X}$.
A straightforward argument using the fact that \bfPplus\ is obtained from \bfP\ by 1-reflecting $X$ downward from $a$ to $b$ (Definition 9.1 of [\ref{refC2}]) shows that $h^+$ is a covering of \bfPplus\ which extends $h$ above $\varphi$.

Suppose \bfPplus\ is obtained from \bfP\ by 2-reflecting $X$ downward from $b$ to $a$.
The proof that $\bfP | \bfPplus$ is cofinally valid in \calB\ is analogous to  the proof of Lemma 9.6 of [\ref{refC2}] and similar to the proof in the previous paragraph (using Definition 9.4 of [\ref{refC2}] instead of Definition 9.1).

Assume $\bfP|\bfPplus$ is a generating rule which is cofinally valid in \calB\ and $\bfP|\bfPstar$ is obtained by 2-reflecting $\bfP|\bfPplus$ upward from $a$ to $b$.
The proof that $\bfP|\bfPstar$ is cofinally valid in \calB\ is analogous to the proof of Lemma 10.3 of [\ref{refC2}].
Let $X=|\bfPplus|\setminus|\bfP|$.
By Definition 10.1 of [\ref{refC2}], \bfPplus\ is a continuous extension of \bfP\ at $a$ (see Definitions 7.1 and 7.4 of [\ref{refC2}]) and $a\preq_2^{\bf P}b$.
Assume $h$ is a covering of \bfP\ in \calB\ and $\varphi$ is a regressive function on the nonminimal indecomposables in the range of $h$.
Since $a\preq_2^{\bf P}b$, $h(a)\preq_2^{\cal B}h(b)$ implying $h(a)\preq_2^\infty h(b)$ in \calB.
Since $\bfP|\bfPplus$ is cofinally valid in \calB, there are cofinally many $\tilde{X}$ below $h(a)$ such that $h[(-\infty,a)^{\bf P}]\cup\tilde{X}$ is closed and a covering of $(-\infty,a)^{\bf P}\cup X$ (as a substructure of \bfPplus).
Since $h(a)\preq_2^\infty h(b)$ in \calB, there are cofinally many $\tilde{X}$ below $h(b)$ such that $h[(-\infty,a)^{\bf P}]\cup \tilde{X}$ is closed and a covering of $(-\infty,a)^{\bf P}\cup X$. Choose such $\tilde{X}$ such that $\varphi(h(b))<\tilde{X}$. 
A straightforward argument using the fact that $\bfP|\bfPstar$ is obtained by 2-reflecting $\bfP |\bfPplus$ upward from $a$ to $b$ (Definition 10.1 of [\ref{refC2}]) shows that $h^+$ is a covering of \bfPstar\ which extends $h$ above $\varphi$.

Assume $\bfPi|\bfPiplusone$ is a generating rule which is cofinally valid in \calB\ for $i<n$  and  \bfPplus\ is a closed substructure of \bfPn\ which extends \bfPzero. 
An easy argument by induction shows $\bfPzero|\bfPi$ is cofinally valid in \calB\ for $i\leq n$. 
The fact that $\bfPzero |\bfPn$ is cofinally valid in \calB\ clearly implies that $\bfPzero|\bfPplus$ is also.

Assume $\bfP|\bfPplus$ is a generating rule which is cofinally valid in \calB\ and $h$ is a continuous embedding of \bfP\ in \bfQ. 
Let \bfQplus\ be a minimal lifting (Definitions 12.1 and 12.4 of [\ref{refC2}]) of $\bfP|\bfPplus$ to \bfQ\ with respect to $h$ and let $h^+$ be the lifting map.
The proof that $\bfQ|\bfQplus$ is cofinally valid in \calB\ is analogous to the proof of Lemma 13.8 of [\ref{refC2}].
By identifying \bfP\ and \bfPplus\ with their images under $h^+$, we may assume $h^+$ is the identity on \unibfPplus.  
Assume $f$ is a covering of \bfQ\ in \calB\ and assume $\varphi$ is a regressive function on the nonminimal indecomposables in the range of $f$.
By increasing the values of $\varphi$ if necessary, we may assume that $f[(-\infty,a)^{\bf Q}]\leq \varphi(h(a))$ whenever $a\in\unibfP$ and $h(a)$ is indecomposable.
Since $\bfP|\bfPplus$ is cofinally valid in \calB, there is a covering $g$ of \bfPplus\ in \calB\ which extends the restriction of $f$ to \unibfP\ above the restriction of $\varphi$ to the indecomposables in $f[\unibfP]$.
The restriction of $f\cup g$ to the idecomposables of \bfQplus\ is order preserving.
By Lemma 4.5 of [\ref{refC2}], this map extends to a unique arithmetic embedding of the arithmetic part of \bfQplus\ in \calB\ which must extend both $f$ and $g$.
Therefore, $f\cup g$ is an arithmetic embedding of the arithmetic part of \bfQplus\ in \calB.
Let \bfQstar\ be the pattern with the same arithmetic part as \bfQplus\ which is induced by \calB\ through $f\cup g$ i.e. so that $f\cup g$ is an embedding of \bfQstar\ in \calB.
Consider the structure \bfQprime\ which has the same arithmetic part as \bfQplus\ so that the interpretation of $\preq_k$ is the intersection of the interpretations of $\preq_k$ in \bfQplus\ and \bfQstar.
A straightforward argument shows \bfQprime\ is a lifting of $\bfP|\bfPplus$ to \bfQ.
Since \bfQplus\ is a minimal lifting of $\bfP|\bfPplus$ to \bfQ,  \bfQprime\ must be a cover of \bfQplus\ (actually, equal to \bfQplus) implying \bfQstar\ is a cover of \bfQplus.
Therefore, $f\cup g$ is a covering of \bfQplus\ in \calB. 
Clearly, $f\cup g$ extends $f$ above $\varphi$.
 \qed

\begin{lem}
\label{l cov ext}
Assume \bfP\ and \bfQ\ are patterns and \bfP\ generates \bfQ.
Any covering of \bfP\ in \calB\ extends to a covering of \bfQ\ in \calB.
\end{lem}
{\bf Proof.} 
Straightforward from the previous lemma (see Definition 14.2 of [\ref{refC2}]).
\qed

\blankline

The following two lemmas will be used only to show that if the arithmetic part of \calB\ is the restriction of \calR\ to an indecomposable of \calR\ then every \calB-covered pattern is covered i.e. \calRtwo-covered. Hence, if one is willing to accept the assumption that every pattern is covered (which increases the proof-theoretic strength of the metatheory to just beyond ${\sf KP}\ell_0$) then these lemmas can be omitted.

The next lemma is an observation that the proofs of parts 3, 4, 6 and 8 of Lemma 14.8 in [\ref{refC2}] actually prove stronger assertions. Notice that in our base theory ${\sf KP}\omega$, saying that a linear ordering is order isomorphic to an ordinal is stronger than saying it is a well-ordering.

\begin{lem} 
\label{l main}
Assume \bfPn\ ($n\in\omega$) is  an increasing sequence of patterns such that $\bfPn | \bfPnplusone$ is a generating rule for each $n\in\omega$. Let $\bfPinfty$ be the union of the \bfPn\ ($n\in\omega$). 
\begin{enumerate}
\item[3*.] 
Every covering of \bfPzero\ in \calB\ extends to a covering of \bfPinfty\ in \calB.
\item[4*.]
Assume  $(|\calB|, \preq^{\cal B})$ is order isomorphic to an ordinal.
If \bfPzero\ is \calB-covered then  $(|\bfPinfty|,\preq^{{\bf P}_\infty})$ is order isomorphic to an ordinal 
\item[6*]
If \bfPinfty\ is a well-ordered structure (i.e. $\preq^{{\bf P}_\infty}$ is a well-ordering of \bfPinfty) and \bfQ\ is a closed substructure of \bfPinfty\ which is a covering of \bfPn\ then $|\bfPn|\preq_{pw}^{{\bf P}_\infty}|\bfQ|$.
\item[8*.]
Assume \bfPn\ ($n\in\omega$) is fair and  \bfPinfty\ is a well-ordered structure.
\begin{enumerate}
\item
For $k=1,2$ and $a,b\in |\calR|$
$$a\preqkinfty b \ \ \ \Longrightarrow \ \ \ a\preq_k^{{\bf P}_\infty}b$$
\item
If $(|\bfPinfty|,\preq^{{\bf P}_\infty})$ is order isomorphic to an ordinal then \bfPinfty\ is isomorphic to $\calRtwo \! \upharpoonright \! \delta$ for some $\delta$ which is indecomposable in \calR.
\end{enumerate}
\end{enumerate}
\end{lem}
{\bf Proof.}
Part 3* follows from Lemma \ref{l cov ext}.

Part 4* follows from part 3*.

For part 6*, notice that parts 1, 5 and 7 of Lemma 14.8 of [\ref{refC2}] implies that \bfPinfty\ satisfies our preliminary assumptions on \calB\ i.e. the arithmetic part of \bfPinfty\ is an arithmetic structure with respect to \calR\ and parts (a)-(c) of the Categoricity Theorem hold. 
Taking \calB\ to be \bfPinfty\ in part 3* we see there is a covering of \bfPinfty\ into itself which extends the covering of \bfPn\ onto \bfQ. Since \bfPinfty\ is well-ordered, we must have $|\bfPn|\preq_{pw}^{{\bf P}_\infty}|\bfQ|$.

The proof of part 8 of Lemma 14.8 of [\ref{refC2}] actually shows part 8*(a) if we replace applications of  part 6 of Lemma 14.8 by applications of part 6* above.

For part 8*(b), we may assume the arithmetic part of \calB\ is $\calR \upharpoonright\delta$ for some ordinal $\delta$ which is indecomposable in \calR\ by parts 1 and 5 of Lemma 14.8 and Lemmas 4.4 and 4.5 of [\ref{refC2}]. 
A simple induction using part 7 of Lemma 14.8 of [\ref{refC2}] and part 8*(a) shows that for $\alpha\leq\delta$, the restriction of $\preq_k^{\cal B}$ to $\alpha$ is the same as the restriction of $\preq_k^{{\cal R}_2}$ to $\alpha$ for $k=1,2$.
\qed

\begin{lem}
If the arithmetic part of \calB\ is isomorphic to an initial segment of \calR\  then any \calB-covered pattern is covered.
\end{lem}
{\bf Proof.}
Assume $h$ is a covering of the pattern \bfP\ in \calB.
Let \bfPn\ ($n\in\omega$) be a fair sequence of patterns with $\bfPzero=\bfP$.

By part 3* of the previous lemma, there is a covering $h^+$ of \bfPinfty\ in \calB\ which extends $h$.
By part 8*(b) of the previous lemma, \bfPinfty\ is isomorphic to an initial segment of \calRtwo.
The restriction of that isomorphism to $|\bfP|$ is a covering of \bfP\ in \calRtwo.
\qed

\blankline

\noindent{\bf Proof of the Categoricity Theorem.}
Our proof will follow the general lines of the proof of Theorem 9.1 of [\ref{refC1}]. 

\halfblankline

{\bf Claim1.} Assume \bfP\ is \calB-covered and \bfPprime\ is a minimal element with respect to $\preq_{pw}^{\cal B}$ (the pointwise partial ordering of finite subsets of \calB) among the closed substructures of \calB\ which are coverings of \bfP.
\begin{enumerate}
\item[(i)]
If \bfQ\ is a substructure of \calB\ which is a cover of \bfP\ then $|\bfPprime| \leqpw |\bfQ|$.
\item[(ii)] 
$\bfP \cong \bfPprime$.
\end{enumerate}

For (i), suppose \bfQ\ is a substructure of \calB\ which is a cover of \bfP. 
By Theorem 14.10 of [\ref{refC2}], there are finite closed substructures $\bf R$ and \bfPstar\ of \calRtwo\ which are isominimal in \calRtwo\ and isomorphic to $\bfPprime\cup\bfQ$ (with a slight abuse of notation) and \bfP\ respectively. 
Let $\overline{\bf P'}$ and $\overline{\bf Q}$ be the images of \bfPprime\ and \bfQ\ respectively under the isomorphism of $\bfPprime\cup\bfQ$ and $\bf R$.
By part 2 of Theorem 14.10 of [\ref{refC2}], $|\bfPstar|\leqpw |\overline{\bf P'}|,|\overline{\bf Q}|$.
By part 5 of Theorem 14.10 of [\ref{refC2}], $\overline{\bf P'}\cup \overline{\bf Q}$ generates $\bfPstar \cup \overline{\bf P'} \cup \overline{\bf Q}$.
By Lemma \ref{l cov ext}, there is a covering $h$ of $\bfPstar \cup \overline{\bf P'} \cup \overline{\bf Q}$ in \calB\ which extends the isomorphism of $\overline{\bf P'}\cup \overline{\bf Q}$ with $\bfPprime\cup\bfQ$. 
Let ${\bf P}''$ be the image of \bfPstar\ under $h$.
We have $|{\bf P}''| \leqpw |\bfPprime|, |\bfQ|$. 
By the minimality of \bfPprime, ${\bf P}'' = \bfPprime$.
Therefore, $|\bfPprime|\leqpw|\bfQ|$.

For part (ii), follow the argument for part (i) (one may take  $\bfQ=\bfPprime$) to conclude from ${\bf P}''=\bfPprime$ that $\bfPstar=\overline{\bfPprime}$. 
Since $\bfPstar\cong\bfP$ and $\overline{\bfPprime}\cong \bfPprime$, $\bfP\cong\bfPprime$.

\halfblankline

For any covered pattern \bfP, let \bfPstar\ be the isominimal substructure of \calRtwo\ which is isomorphic to \bfP.
For \bfP\ an isominimal substructure of \calB, define $f_{\bf P}$ to be the isomorphism of \bfP\ and \bfPstar. 
Let $f$ be the union of the $f_{\bf P}$. 
A straightforward argument shows $f$ is an embedding of the core of \calB\ into the core of \calRtwo.

To show the range of $f$ is an initial segment of \calRtwo, assume $\alpha<\beta$ where $\beta$ is in the range of $f$.
There is an isominimal substructure \bfP\ of \calB\ such that $\beta$ is in the range of $f_{\bf P}$.
Let \bfPn\ ($n\in\omega$) be a fair sequence with $\bfPzero=\bfP$ and let \bfPinfty\ be the union of the \bfPn.
By Lemma 14.9 of [\ref{refC2}], there is an isomophism $g$ of \bfPinfty\ with $\calRtwo \! \upharpoonright \! \delta$ for some $\delta$ which is indecomposable in \calR\ and the image of \bfPn\ under $g$ is \bfPnstar\ for each $n\in\omega$.
Fix $n$ such that $\alpha$ is in \bfPnstar.
By Lemma \ref{l cov ext} and Claim 1, there is an isominimal substructure \bfQ\ of \calB\ which is isomorphic to \bfPn.
Since $\alpha$ is in \bfPnstar\ which  the range of $f_{\bf Q}$, $\alpha$ is in the range of $f$.
\qed

\begin{cor}
\label{c characterization}
Assume \calB\ satisfies (a)-(d) of the Categoricity Theorem for \calRtwo.
If there are arbitrarily long finite chains in $\preq_2^{\cal B}$  then the core of \calB\ is isomorphic to the core of \calRtwo\ and a finite structure is isomorphic to a finite closed substructure of \calB\  iff it is a pattern of resemblance of order 2. 
\end{cor}
{\bf Proof.}
Assume there are arbitrarily long finite chains in $\preq_2^{\cal B}$. 
By the Categoricity Theorem, the core of \calB\ is isomorphic to an initial segment of the core of \calRtwo.
Since this initial segment contains arbitrarily long finite chains in \leqtwo, it must be the entire core of \calRtwo\ by part 2 of Theorem 14.10 of [\ref{refC2}].
Hence, every pattern of resemblance of order two is isomorphic to a substructure of \calB.
The converse is straightforward after noticing that condition (a) of the Categoricity Theorem implies that $\alpha$ is indecomposable whenever $\alpha \preq_1^{\cal B} \beta$ and both $\alpha$ and $\beta$ are indecomposable whenever $\alpha \preq_2^{\cal B} \beta$. 
\qed

\begin{cor}
Assume ${\cal R}_2'$ is the alternate definition of \calRtwo\ from the introduction to {\rm [\ref{refC2}]} using $\Sigma_1$ and $\Sigma_2$ elementarity.
${\cal R}_2'\! \upharpoonright \! \delta$ satisfies the (a)-(d) of the Categoricity Theorem for each indecomposable $\delta$ and, hence, the conclusions of the Categoricity Theorem and the previous corollary hold for \calRtwo.
\end{cor}
{\bf Proof.} Straightforward after noting that in ${\cal R}_2'$, if $\alpha<\beta$, $\beta$ is a limit ordinal and $\alpha \leqone \xi$ for all $\xi$ with $\alpha \leq \xi < \beta$ then $\alpha\leqone \beta$.
\qed

\blankline

One can prove that \leqtwo\ in ${\cal R}_2'$ has arbitrarily long finite chains well within $\sf ZF$.

\blankline

\begin{center}
{\bf References}
\end{center}
\begin{enumerate}
\item 
\label{refC1}
{\it Elementary patterns of resemblance,} Annals of Pure and Applied Logic 108 (2001), pp. 19-77.
\item
\label{refC2}
{\it Patterns of resemblance of order 2,} Annals of Pure and Applied Logic 158 (2009), pp. 90-124.
\end{enumerate}
\end{document}